\documentclass[11pt]{amsart}

\usepackage{amssymb}

\newtheorem{theorem}{{Theorem}}[section]
\newtheorem{defn}[theorem]{{Definition}}%[section]
\newtheorem{proposition}[theorem]{{Proposition}}%[section]
\newtheorem{lemma}[theorem]{{Lemma}}%[section]
\newtheorem{corollary}[theorem]{{Corollary}}%[section]
\newtheorem{propn}[theorem]{{Proposition}}
\newtheorem{fact}[theorem]{{Fact}}%[section]
\newtheorem{remark}[theorem]{{Remark}}%[section]
\theoremstyle{definition}
\newtheorem{example}[theorem]{{Example}}

\newcommand{\ad}{\overline{\mathrm{ad}}}
\newcommand{\R}{{\bf R}}
\newcommand{\Z}{{\bf Z}}

\begin{document}

\title{Actions of noncompact semisimple groups on Lorentz manifolds}
\author{M. Deffaf, K. Melnick  and A. Zeghib}
\date{\today}
\maketitle

\begin{abstract}

The above title is the same, but with ``semisimple'' instead of ``simple,''  as that of 
 a notice by Nadine Kowalsky. There, she announced many theorems on the subject of actions of simple Lie groups preserving
a Lorentz structure. Unfortunately, she published proofs for essentially only half of the announced results before her premature death. 
Here, using a different, geometric approach, we generalize her results to the semisimple case, and give proofs of all her announced results.
\end{abstract}

\section{Introduction}

Isometric actions on Lorentz manifolds were first investigated in the 
compact case (see \cite{ZiLor}, \cite{Gr}, \cite{AS1}, \cite{AS2}, \cite{Ze3}, \cite{Ze4}).  The natural question was then: how can a compact Lorentz
manifold have a noncompact  isometry group?
 There is strong evidence that such a question is in fact ``decidable'' for a wide class of geometric structures (see, for instance \cite{DAG}, \cite{Zirigid}).

\subsection{Framework} One new aspect of Kowalsky's work was to deal with actions of groups on {\it noncompact}
Lorentz manifolds. 
Obviously, nothing can be said about such actions
 without compensating for noncompactness with a dynamical counterpart ensuring some kind of recurrence.  A natural and rather weak condition used by Kowalsky is nonproperness of the action. 

Let us here appreciate consideration of the noncompact case, at least from 
a physical point of view, according to which compact spacetimes have little interest.
Having a nonproper isometry group is a manifestation of a non-Riemannian character of the geometry of 
spacetime. It is in such spaces that one can observe ``dilation of length'' and ``contraction of time.'' It is surely interesting to try to classify spacetimes with nonproper isometry groups.
This job, however, does not seem to be easy. Some extra hypotheses are therefore in order.  Kowalsky restricts to actions of {\it simple Lie groups}.

\subsection{Kowalsky's main theorem} 

The de Sitter and anti-de Sitter spaces, $dS_n$ and $AdS_{n+1}$, respectively, are the homogeneous spaces $O(1, n)/ O(1, n-1)$ and $O(2, n)/O(1, n)$. Geometrically, they are the universal Lorentz spaces of constant positive and negative curvature, respectively. A striking fact proved by Kowalsky is that,  {\it at a  group level}, they are the only Lorentz nonproper $G$-spaces, with $G$ simple:

\begin{theorem} [Kowalsky \cite{Kow1} 5.1]
 \label{main.kowalsky}
Let $G$ be a simple Lie group with finite center acting isometrically and nonproperly on a connected Lorentz manifold.  Then $G$ is locally isomorphic to either $O(1,n), n \geq 2,$ or $O(2,n), n \geq 3$.
\end{theorem}

\begin{remark} 
The groups $O(1,1)$ and  $O(2,2)$ are not simple.
\end{remark}

\subsection{Geometry of semisimple isometric actions} 

Once the acting group is known, the problem arises to understand the geometry of the Lorentz space, or at least that of orbits.  Here, one hopes the space looks like de Sitter or anti-de Sitter space, depending on whether $G$ is locally $O(1,n)$ or $O(2,n)$.  Nadine Kowalsky announced results to this effect in \cite{Kow2} and wrote proofs for the $O(1,n)$ case in her thesis \cite{Kowthesis}.  We will recall their statements below in \S \ref{legacy}.  Unfortunately, she prematurely died, before publishing proofs.  

\subsection{The technique} When a Lie group $G$ acts on $M$ preserving a 
pseudo-Riemannian metric, one can consider a Gauss map from $M$ to $S^2(\mathfrak{g})$, the space of quadratic forms on the Lie algebra $\mathfrak{g}$ of $G$. When $S^2(\mathfrak{g})$ is endowed with the natural adjoint 
$G$-action, the Gauss map is equivariant, and the $G$-space $S^2(\mathfrak{g})$ reflects the dynamics on $M$. It is via this map that the non-properness condition is translated as a geometric condition on the induced metrics on orbits. This idea, due to Kowalsky, has become a basic tool in similar questions on the subject, e.g., Adams-Stuck \cite{AS1}, \cite{AS2}, and Bader-Nevo \cite{BN} (Remark here that variants of the Gauss map, with other natural spaces instead of 
$S^2(\mathfrak{g})$ were used by other authors, e.g. Gromov \cite{Gr} and Zimmer \cite{Zirigid}).
However, this is the starting point; further work in the proof is algebraic Lie theoretic. 

\subsection{Other works} Another proof of Theorem \ref{main.kowalsky} was proposed by S. Adams \cite{Adams}; his methods involve an analysis similar to Kowalsky's, except that zero jets of germs of Killing vector fields are replaced by higher order jets.  In other directions,
  S. Adams investigated  Lorentz-isometric actions, for simply connected Lie groups, with the stronger dynamical condition that some orbit either is not closed or has noncompact stabilizers (\cite{Adamsonp1}, \cite{Adamsonp2}).
  
  Concerning Kowalsky's unpublished proofs, we notice a contribution by D. Witte Morris \cite{wm}, in which he considers the homogeneous case.  More precisely, he takes $G$ locally isometric to the isometry group of de Sitter or anti-de Sitter space, respectively, and considers a nonproper Lorentz homogeneous space $G/H$---that is, $H$ is noncompact, and the $G$-action on $G/H$ preserves some Lorentz metric. He proves that $\mathfrak{h}$ is isomorphic to a standard copy of $\mathfrak{o}(1,n-1)$ in $\mathfrak{o}(1,n)$ or of $\mathfrak{o}(1,n)$ in $\mathfrak{o}(2,n)$; it follows that $G/H$ is locally homothetic to de Sitter or anti-de Sitter space.  Witte Morris' proof is quite algebraic.
 
 \subsection{On the present contribution} \label{contribution}
 Our investigation here highly relies on the approach of \cite{ADZ}, although the two articles can be read completely independently.
 From \cite{ADZ}, we will use the result recalled below as Theorem \ref{adz}.  In light of this result, on the structure of nonproper orbits of Lorentz type, the present paper addresses the case in which the acting group has a nonproper degenerate orbit.

 Before stating our results, let us give some motivations and emphasize new features:
 
  $\bullet$ {\it Completing Kowalsky's work:} One major goal is to prove the announced results of Nadine Kowalsky. 
 
  $\bullet$ {\it Geometric approach: } The approach here is different from
  that of Kowalsky (as well as from others', for instance Adams'). Together with \cite{ADZ}, we get proofs of the main results, in particular, of Theorem \ref{main.kowalsky}, using many geometric arguments, where one sees the global structure of proofs.

   $\bullet$ {\it From simple to semisimple:}  More important, we generalize results to the semisimple case, assuming there are no local $SL_2(\R)$-factors. A semisimple Lie group is essentially a product of simple Lie groups, but, in general, a nonproper action of a product does not derive from a nonproper action of one factor.  However, in the Lorentz setting, we conclude that it necessarily does---that is, the semisimple case reduces to the simple one.  This is really an important fact, since it leads one to hope to reduce the remaining work to the case in which the group is solvable.  Of course, the reason for this is the Levi decomposition of Lie groups, which says that a Lie group is essentially a semidirect product of a semisimple and a solvable group. 
   
\subsection{Kowalsky's legacy} \label{legacy} 

In \cite{Kow2}, the following theorems are stated.  For Theorem \ref{Kow.stabs1} below, see also \cite{Kowthesis} 6.2.  The manifold $M$ and group $G$ are assumed connected throughout.

\begin{theorem}
\label{Kow.stabs1}
Let $G$ be locally isomorphic to $O(1,n), n \geq 3$, and suppose that $G$ acts on a manifold $M$ preserving a Lorentz metric.  Then all noncompact stabilizers have a Lie algebra isomorphic to either $\mathfrak{o}(1,n)$, $\mathfrak{o}(1,n-1)$, or $\mathfrak{o}(n-1) \ltimes \R^{n-1}$.
\end{theorem}

\begin{theorem}
\label{Kow.stabs2}
Let $G$ be locally isomorphic to $O(2,n), n \geq 6$, with $G$ having finite center.  Suppose that $G$ acts nontrivially on a manifold $M$ preserving a Lorentz metric.  Then all noncompact stabilizers have a Lie algebra isomorphic to $\mathfrak{o}(1,n)$.
\end{theorem}

\begin{theorem}
\label{Kow.warped.product}
Let $G$ and $M$ be as in Theorem \ref{Kow.stabs2} above, and assume there is a point of $M$ with noncompact stabilizer.  Then the universal cover $\widetilde{M}$ is Lorentz isometric to a warped product $L \times_w \widetilde{AdS}_{n+1}$, where $\widetilde{AdS}_{n+1}$ is the simply connected $(n+1)$-dimensional Lorentz space of constant curvature $-1$, and $L$ is a Riemannian manifold.  Further, the induced action of the universal cover $\widetilde{G}$ on $\widetilde{M}$ is via the canonical $\widetilde{G}$-action on $\widetilde{AdS}_{n+1}$ and the trivial action on $L$.
\end{theorem}

See Section \ref{background} below for the definition of warped product.

 \subsection{Results} As said above in \S \ref{contribution}, we provide here proofs of all previous statements of Kowalsky, together with some improvements. 

A submanifold $N$ in a Lorentz manifold is \emph{degenerate} if $T_xN^{\perp} \cap T_xN \neq {\bf 0}$.  In Minkowski space $\R^{1,n}$, the simple subgroup $O(1,n) \subset Isom(\R^{1,n})$ has one degenerate orbit, which, together with the origin, forms the \emph{light cone}, the set of all isotropic vectors in $\R^{1,n}$.  The stabilizer in $O(1,n)$ of a nonzero vector in the light cone is isomorphic to $O(n-1) \ltimes \R^{n-1}$, where the action of $O(n-1)$ on $\R^{n-1}$ is the usual representation.

In the degenerate case, we have the following theorem, which says that a degenerate orbit for a simple group acting nonproperly is locally homothetic to the Minkowski light cone.  Together with Theorem 1.5 of \cite{ADZ}, which classifies nonproper orbits of Lorentz type, it implies Theorems \ref{Kow.stabs1} and \ref{Kow.stabs2} of Kowalsky above.

\begin{theorem}
\label{degstabs}
Suppose $G$ is a connected group with finite center, locally isomorphic to $O(1,n)$ or $O(2,n)$ for $n \geq 3$.  If $G$ acts isometrically on a Lorentz manifold and has a degenerate orbit $O$ with noncompact stabilizer $G(x)$, then 
\begin{enumerate}
\item{$G$ is locally isomorphic to $O(1,n)$}
\item{The Lie algebra $\mathfrak{g}(x)$ is isomorphic to $\mathfrak{o}(n-1) \ltimes \R^{n-1}$.} 
\item{The orbit $O$ is locally homothetic to the light cone in Minkowski space.}
\end{enumerate} 
\end{theorem}

The following result implies Theorem \ref{Kow.warped.product} above.

\begin{theorem} \label{warped}

   If $G$, a group with finite center locally isomorphic to  $O(2, n)$, $n \geq 3$,
   acts isometrically and nontrivially on a Lorentz manifold $M$, with some noncompact stabilizer, then, up to finite covers, 
   $M$ is equivariantly isometric to a warped product $L\times_w AdS_{n+1}$ of a Riemannian
   manifold $L$ with the anti-de Sitter space $AdS_{n+1}$.
 \end{theorem}

We extend the above results to semisimple groups. Note that the noncompact stabilizer assumption is weakened to nonproperness of the action. A local factor of a semisimple Lie group $G$ is a Lie group $G_1$ such that the Lie algebra $\mathfrak{g}_1$ is a direct summand of $\mathfrak{g}$.  If $|Z(G)| < \infty$ and $G$ acts on a manifold $M$, then a finite cover of any local factor $G_1$ acts on $M$.  We will slightly abuse terminology below by referring to the ``$G_1$-action'' and ``$G_1$-orbits'' on $M$.

\begin{theorem} \label {main.result} Let $G$ be a semisimple group with finite center and no local $SL_2(\R)$-factor, acting isometrically, faithfully, and nonproperly on a Lorentz manifold $M$. Then 

\begin{enumerate}
\item{$G$ has a local factor $G_1$ isomorphic to $O(1, n)$ or $O(2, n)$.}
\item{There exists a Lorentz manifold $S$, isometric, up to finite cover, to $dS_n$ or $AdS_{n+1}$, depending whether $G_1$ is isomorphic to $O(1,n)$ or $O(2,n)$, and an open subset of $M$, in which each $G_1$-orbit is homothetic to $S$.  }
\item{Any such orbit as above has a $G_1$-invariant neighborhood isometric to a warped product $L \times_w S$, for $L$ a Riemannian manifold.}
\end{enumerate}
\end{theorem}

For $O(1,n)$-actions, we can also describe orbits with compact isotropy (Proposition \ref{partition}).  Nonproper $O(1,n)$-actions for which compact isotropy occurs strongly resemble the standard action on ${\bf R}^{1,n}$.  As a fusion, we can give the following ``full'' theorem:

\begin{theorem} \label{full}  Let $G$ be a semisimple Lie group with finite center and no local $SL_2(\R)$-factor, acting nonproperly and isometrically on a Lorentz manifold $M$.  Then, $G$ has a simple local factor $G_1$ that acts nonproperly.  There are two possibilities for $G_1$:

\begin{enumerate}
\item{$G_1 \cong O(2, n)$. In this case, there is a Lorentz manifold $S$ isometric, up to finite cover, to $AdS_{n+1}$, such that all $G_1$-orbits are homothetic to $S$.
In fact, up to finite cover, $M$ is a warped product $L \times_w AdS_{n+1}$.}

\item{$G_1 \cong O(1, n)$. There are open sets $U$ and $V$ such that $M = U \sqcup \partial U \sqcup V$, where }

\begin{itemize}

\item{For any $x \in U$, there exists $S$, isometric to $dS_n$ up to finite cover, such that the $G_1$-orbit of $x$ is homothetic to $S$. A neighborhood of $x$ is $G_1$-equivariantly isometric to a warped product $L \times_w S$, for some Riemannian manifold $L$.}

\item{Orbits of $G_1$ on the boundary of $U$ are either fixed points or locally homothetic to the light cone in Minkowski space $\R^{1, n}$; further $\partial U = \partial V$.}

\item {For any $x \in V$, the $G_1$-orbit of $x$ is homothetic to ${\bf H}^n$.  The set $V$ is globally $G_1$-equivariantly isometric to a warped product $L \times_w {\bf H}^n$, for some Lorentz manifold $L$. }

\end{itemize}
\end{enumerate}

\end{theorem}

The following two examples of nonproper $SL_2(\R)$-actions illustrate the necessity of the hypothesis of no $SL_2(\R)$-factors.

\begin{example} 
\label{ex1sl2}
 In this example, the manifold has constant curvature, but has no common finite cover with any of the constant-curvature models.  The group $PSL_2(\R)$ with the Lorentz metric arising from the Killing form is isometric to $AdS_3$.  For $\Gamma$ a cocompact lattice in $PSL_2(\R)$, the manifold $M = AdS_3 / \Gamma$ admits an isometric, nonproper left-action by $PSL_2(\R)$.  
\end{example}

\begin{example} 
\label{ex2sl2}
 This example is a transitive, nonproper, isometric $SL_2(\R)$-action on a manifold with nonconstant curvature.  Let $g$ be the Killing form on $\mathfrak{sl}_2(\R)$.  Let $\lambda_1, \lambda_2$ and $\lambda_3$ be the standard basis for $(\R^3)^*$.  There is an isomorphism $\mathfrak{sl}_2(\R) \cong \R^3$ for which the Killing form is $2 \lambda_1 \cdot \lambda_2 + \lambda_3^2$.  Under this isomorphism, an $\R$-split element of $\mathfrak{sl}_2(\R)$ maps to the basis element ${\bf e}_3$.  Let $c \neq 1$ be a positive number.  Let $g^{\prime}$ be the pullback of the form $2 \lambda_1 \cdot \lambda_2 + c \lambda_3^2$ to $\mathfrak{sl}_2(\R)$.  The adjoint action of the $\R$-split torus $A$ of $SL_2(\R)$ preserves $g^{\prime}$.  Let $\Gamma$ be a lattice in $A$.  The form $g^{\prime}$ gives rise to an $SL_2(\R)$-invariant Lorentz metric on $M = SL_2(\R) / \Gamma$.  The isometric $SL_2(\R)$-action is nonproper, but $M$ does not have constant curvature.
\end{example}

For Lorentz-isometric actions of $SL_2(\R)$ on finite volume manifolds, Gromov has shown that all stabilizers are discrete, and the universal cover is isometric to a warped product $L \times_w \widetilde{AdS}_3$ (\cite{Gr} 5.4.A).

This final example illustrates the necessity for the results above of the hypothesis that $G$ has finite center.  We thank the referee for bringing this example to our attention.
\begin{example}
\label{ex.infinite.center}
The main result of \cite{AdInduction} implies that, if $Z$ is a central closed subgroup of a Lie group $G$, and if $Z$ acts isometrically on a Lorentz manifold $M$, then the manifold $G \times_Z M$ admits a $G$-invariant Lorentz metric.  There exists a Lorentz manifold $M$, diffeomorphic to a torus, with an Anosov action of the integers $\Z$ by isometries (see \cite{BarbZe} \S 6.3.3); in fact, this action has a fixed point.  Then $\widetilde{O}(2,n)$ acts isometrically and nonproperly on $\widetilde{O}(2,n) \times_{\Z} M$.  All stabilizers are discrete, and there is obviously no warped product as in the theorems above. 
\end{example}

\emph{Notation and Terminology}: The Lie algebra of a Lie group $G$ will be denoted $\mathfrak{g}$.  The stabilizer subgroup in $G$ and corresponding subalgebra of a point $x$ will be denoted $G(x)$ and $\mathfrak{g}(x)$, respectively.

Lie groups $G$ and $H$ are \emph{locally isomorphic} if $\mathfrak{g} \cong \mathfrak{h}$.  As discussed above, a group $G_1$ is a \emph{local factor} of $G$ if $\mathfrak{g}_1$ is a direct summand of $\mathfrak{g}$.

The dimension of the Lorentz manifold $M$ will be denoted $d$ throughout.

\section{Background: warped product near Lorentz orbits} \label{background}

\begin{defn}
For two pseudo-Riemannian manifolds $(L, \lambda)$ and $(S,\sigma)$, a \emph{warped product} $L \times_w S$ is given by a positive function $w$ on $L$ : the metric at $(l,s)$ is $\lambda_l + w(l) \sigma_s$. 
\end{defn}

\subsection{Results of \cite{ADZ}}   We will make use of the following theorem:
 
\begin{theorem} \label{adz} (\cite{ADZ} 1.5)
Let  $G$ be a connected semisimple Lie group acting isometrically on a Lorentz manifold $M$ of dimension $\geq 3$.  Suppose  that no local factor of $G$ is isomorphic to $SL_2(\R)$ and that there exists an orbit $O$ of Lorentz type with noncompact isotropy.

Then, up to a finite cover, $G$ factors $G \cong G_2 \times G_1$, where:

\begin{enumerate}
\item{$G_1$ possesses an orbit $O_1$ which is a Lorentz space of constant, non-vanishing curvature, and 
$G_1$ equals $Isom^0(O_1)$.}

\item{There is a $G$-invariant neighborhood $U$ of $O_1$ which is a warped product $L \times_w O_1$.}

\item{The factor $O_1$ corresponds to $G_1$-orbits, and $G_2$ acts along the $L$-factor.}
\end{enumerate}
\end{theorem}

\section{Properties of the isotropy representation}
\label{section.isotropy.props}
Here we collect some algebraic facts about the structure of nonproper degenerate orbits.  Suppose that $x$ is a point of $M$ with degenerate $G$-orbit.  Denote this orbit by $O$.  Recall that $d$ is the dimension of $M$, and assume that $G$ is semisimple.

Fix an isometric isomorphism of $T_x M$ with $\R^{1,d-1}$, determining an isomorphism $O(T_xM) \cong O(1,d-1)$.  Let $V$ be the image of $T_xO$ in $\R^{1,d-1}$.  Let $\Phi : G(x) \rightarrow O(1,d-1)$ be the resulting isotropy representation.  Because $G$ acts properly and freely on the bundle of Lorentz frames of $M$, the isotropy representation is an injective, proper map.  The invariant subspace $V$ is degenerate, so $\Phi(G(x))$ leaves invariant the line $V^{\perp} \cap V$.  The stabilizer of an isotropic line is conjugate in $O(1,d-1)$ to the parabolic 
$$ P =  (K \times A) \ltimes U $$

where $U \cong {\bf R}^{d-2}$ is unipotent, $A \cong {\bf R}^*$, and $K \cong O(d-2)$ with the conjugation action of $K \times A$ on $U$ equivalent to the standard representation of the conformal group of ${\bf R}^{d-2}$.  Denote by $\mathfrak{p}$ the Lie algebra of $P$, and by $\mathfrak{k}$, $\mathfrak{a}$, $\mathfrak{u}$, the subalgebras corresponding to $K$, $A$ and $U$.  Let  $\varphi : \mathfrak{g}(x) \rightarrow \mathfrak{o}(1,d-1)$ be the Lie algebra representation tangent to $\Phi$.

Note that $\mathfrak{g}/ \mathfrak{g}(x)$ can be identified with $T_xO$ by the map 
$$Y \mapsto \overline{Y}(x) = \left. \frac{\partial}{\partial t} \right|_0 e^{tY} x$$

For $g \in G(x)$, differentiating the relation $ge^{tY}x = (ge^{tY}g^{-1})x$ gives 

$$D_x g (\overline{Y}(x)) =  \overline{\mbox{Ad}g(Y)} (x)$$

In other words, $\Phi$ restricted to $V$ is equivalent to the representation $\overline{\mbox{Ad}}$ of $G(x)$ on $\mathfrak{g} / \mathfrak{g}(x)$ arising from the adjoint representation.  Let $\ad$ be the representation tangent to $\overline{\mbox{Ad}}$.  

An element $Y$ of $\mathfrak{g}$ is called \emph{nilpotent} if ad($Y$) is nilpotent.  An element $Y$ is \emph{$\R$-split} if ad$(Y)$ is diagonalizable over $\R$.

\begin{lemma}
\label{3props}
The stabilizer subalgebra $\mathfrak{g}(x) \subset \mathfrak{g}$ has the following properties:

\begin{enumerate}

\item{For all $Y \in \mathfrak{g}(x)$, the endomorphism $\ad (Y)$ has no real nonzero eigenvalues.  In fact, the same is true for $\varphi$, so $\varphi(\mathfrak{g}(x))$ is conjugate to a subalgebra of $\mathfrak{k} \ltimes \mathfrak{u}$.} 

\item{The stabilizer subalgebra $\mathfrak{g}(x)$ contains no element $\R$-split in $\mathfrak{g}$.}

\item{There exists a subalgebra $\mathfrak{s}(x)$ in which $\mathfrak{g}(x)$ has codimension one such that $[\mathfrak{s}(x),\mathfrak{s}(x)] \subset \mathfrak{g}(x)$, and the representation of $\mathfrak{g}(x)$ on $\mathfrak{g}/ \mathfrak{s}(x)$ is skew-symmetric---that is, every endomorphism is diagonalizable with purely imaginary eigenvalues. }  
\end{enumerate} 
\end{lemma}

\begin{proof} 
\ 
\begin{enumerate}

\item{Suppose that $\ad (Y)$ has eigenvalue $\lambda > 0$.  Then $\lambda$ is also an eigenvalue of $\varphi(Y)$ on $V$.  Since $\varphi(Y)$ is skew-symmetric on $V / (V \cap V^{\perp})$, the generalized eigenspace for $\lambda$ in $V$ is one-dimensional and equals $V^{\perp} \cap V$.  The trace of $\varphi(Y)$ on $V$ is $\lambda$, so the trace of $\ad(Y)$ on $\mathfrak{g} / \mathfrak{g}(x)$ is $\lambda$.  

Next we will show that the trace of $\mbox{ad}(\varphi(Y))$ on $\varphi(\mathfrak{g}(x))$ is nonnegative.  Recall that $\mathfrak{p}$ decomposes $(\mathfrak{a} \times \mathfrak{m}) \ltimes \mathfrak{u}$.  The restriction of $\mbox{ad}(\mathfrak{p})$ to $\mathfrak{u}$ factors through the projection to $\mathfrak{a} \times \mathfrak{m}$, which acts as the standard conformal representation on $\R^{d-2}$.  If $\varphi(Y)$ has eigenvalue $\lambda$ on $V \cap V^{\perp}$, then it projects to $\lambda I + M$ in $\mathfrak{a} + \mathfrak{m} \cong \mathfrak{conf}(\R^{d-2})$.  Therefore, the characteristic roots of $\mbox{ad}(\varphi(Y))$ on $\mathfrak{u}$ are in $\lambda + i \R$.  On $\mathfrak{p} / \mathfrak{u} \cong \mathfrak{a} \oplus \mathfrak{m}$, the characteristic roots of any adjoint endomorphism are all purely imaginary.  Therefore, the characteristic roots of $\mbox{ad}(\varphi(Y))$ on $\mathfrak{p}$ are in $(\lambda + i \R) \cup i \R$, implying that the trace of $\mbox{ad}(\varphi(Y))$ is nonnegative on $\varphi(\mathfrak{g}(x))$.

Then the trace of ad$(Y)$ on $\mathfrak{g}(x)$ is also nonnegative.  Finally, the trace of ad$(Y)$ on $\mathfrak{g}$ is positive, contradicting the unimodularity of $\mathfrak{g}$.  If $\lambda < 0$, the same argument shows that the trace of ad$(Y)$ on $\mathfrak{g}$ is negative.  Therefore, no $\overline{\mathrm{ad}}(Y)$ has any nonzero real eigenvalues, and no $\varphi(Y)$ has any nonzero real eigenvalues on $V$.

If $\varphi(Y)$ has a nonzero real eigenvalue on $\R^{1,d-1}$, then an eigenvector must be isotropic.  It either lies in $V$ or is not orthogonal to $V^{\perp} \cap V$.  In either case, $\varphi(Y)$ has a nonzero real eigenvalue on $V$, a contradiction. }

\item{If an $\R$-split element $H \in \mathfrak{g}(x)$, then by (1), all root vectors on which ad$(H)$ is nontrivial must project to ${\bf 0}$ in $\mathfrak{g} / \mathfrak{g}(x)$.  In this case, $\mathfrak{g}(x)$ contains a subalgebra isomorphic to $\mathfrak{sl}_2(\R)$.  Applying the monomorphism $\varphi$ would yield a subalgebra isomorphic to $\mathfrak{sl}_2(\R)$ in $\mathfrak{p}$, which is impossible.}

\item{Take any  $Z^{\prime}$ spanning $V^{\perp} \cap V$.  Item (1) implies $\varphi(\mathfrak{g}(x))$ annihilates $Z^{\prime}$.  Take the corresponding vector in $\mathfrak{g} / \mathfrak{g}(x)$, and let $Z$ be any lift to $\mathfrak{g}$.  Then $\mbox{ad}(\mathfrak{g}(x))(Z) \subseteq \mathfrak{g}(x)$, so $\mathfrak{s}(x) = \R Z + \mathfrak{g}(x)$ is the desired subalgebra.

From the equivalence of $\varphi|_V$ with $\ad$, the representations  $V/(V \cap V^{\perp})$ and $\mathfrak{g}/\mathfrak{s}(x)$ are equivalent.  The former is skew-symmetric.}
\end{enumerate}
\end{proof}

\section{Root spaces in isotropy subalgebra}
\label{section.roots.in.isotropy}

By a nonproper orbit we will mean one with noncompact isotropy.  Theorem 1.8 of \cite{ADZ} asserts the existence of a nonproper orbit under the assumptions of our Theorem \ref{main.result} above.  The proof in \cite{ADZ} was easily deduced from the following result of \cite{Kow1}:

\begin{proposition}  \label{root.isotropic} If the $G$-action is nonproper, then there is $x \in M$, and an $\R$-split element $H$ of $\mathfrak{g}$ such that the negative root space 
$$\Sigma_{\alpha(H) <0} \mathfrak{g}_{\alpha}$$ 
is isotropic at $x$.  If $G$ has noncompact stabilizer at some point $y$ of $M$, then we may take $y = x$ above.
\end{proposition} 

The subalgebra $\mathfrak{s}(x)$ of the previous section is exactly the maximal subspace $\mathfrak{s} \subseteq \mathfrak{g}$ such that $\{ \overline{Y}(x) \in T_xM \ : \ Y \in \mathfrak{s} \}$ is an isotropic subspace.  

Fact 4.2 of \cite{ADZ} is that 
$$\left( \Sigma_{\alpha(H) <0} \mathfrak{g}_{\alpha} \right) \cap \mathfrak{g}(x) \neq {\bf 0}$$ 

If this intersection were ${\bf 0}$, then the subalgebra $\oplus_{\alpha(H) \geq 0} \mathfrak{g}_{\alpha}$ would have codimension one in $\mathfrak{g}$.  Such a subalgebra could only exist if $\mathfrak{sl}_2(\R)$ were a factor of $\mathfrak{g}$, but our hypotheses exclude this.  

\begin{remark}
Point stabilizers are discrete for the $SL_2(\R)$-actions given in Examples \ref{ex1sl2} and \ref{ex2sl2}. 
\end{remark}

Denote by $\mathfrak{a}$ and $\Delta$ the Cartan subalgebra and root system, respectively, of Proposition \ref{root.isotropic}.  The remainder of this section is devoted to showing the following proposition.

\begin{proposition} 
\label{propJ}
There exist $J \in \mathfrak{a}$ and $S \subset \Delta$ such that 

\begin{enumerate}
\item{ $\mathfrak{s}(x) = \R J + \mathfrak{g}(x)$}
\item{ $\alpha(J) < 0$ for all $\alpha \in S$ }
\item{ $\Sigma_{\alpha \in S} \mathfrak{g}_{\alpha} \subseteq \mathfrak{g}(x)$}
\item{ $\mbox{dim}( \Sigma_{\alpha \in S} \mathfrak{g}_{\alpha}) \geq 2$}
\end{enumerate}
\end{proposition}

\begin{proof}

Fix $H$ as in Proposition \ref{root.isotropic}.  Let 
$${\bf 0} \neq X \in \left( \Sigma_{\alpha(H) < 0} \mathfrak{g}_{\alpha} \right) \cap \mathfrak{g}(x) $$  

There exist $J \in \mathfrak{a}$ and a nilpotent $Y$ in $\mathfrak{g}$ such that $[J,X] = - 2X$ and $[X,Y] = J$ (see \cite{sam} 2.4.B).  The operator $\ad(X)$ is nilpotent; on the other hand, by Lemma \ref{3props} (3), $\ad(X)$ is skew-symmetric on $\mathfrak{g} / \mathfrak{s}(x)$, so $\mbox{ad}(X)(\mathfrak{g}) \subseteq \mathfrak{s}(x)$.  Therefore, $J$ belongs to $\mathfrak{s}(x)$.  
 
Note $J \notin \mathfrak{g}(x)$ by Lemma \ref{3props} (2).
 Therefore $\mathfrak{s}(x) = \R J + \mathfrak{g}(x)$, proving (1).

Let $S$ be the set of $\alpha \in \Delta$ such that $\alpha(H) < 0$ and $\alpha(J) < 0$, so (2) is obviously satisfied.  From the relation $[J,X] = -2X$, any $\alpha$ such that $X$ has a nontrivial component in $\mathfrak{g}_{\alpha}$ satisfies $\alpha(J) = -2$, so any such $\alpha$ belongs to $S$; in particular $S$ is not empty.

For $\alpha \in S$, we have $\mathfrak{g}_{\alpha} \subset \mathfrak{s}(x)$ and 
$$\mathfrak{g}_{\alpha} = [J,\mathfrak{g}_{\alpha}] \subset [\mathfrak{s}(x),\mathfrak{s}(x)] \subset \mathfrak{g}(x)$$

by Lemma \ref{3props} (3), showing statement (3) above. 

Now, replacing $X$ by a nonzero element of some $\mathfrak{g}_{\alpha}$, $\alpha \in S$, we may assume that $ - J$ is a \emph{basic translation}---that is, there exists $c_{\alpha} < 0$ such that 
$$\alpha(J) = - 2 \qquad \mbox{and} \qquad \alpha(Z) = c_{\alpha} \kappa(J,Z)$$ 
for any $Z \in \mathfrak{a}$, where $\kappa$ denotes the Killing form.  In this case, we have that for any root $\beta$, the reflection
$$ \sigma_{\alpha}(\beta) = \beta + \beta(J) \alpha$$
is again a root (see \cite{sam} II.5.A).

Now, to show (4) it suffices to show that $\mbox{dim}(\mathfrak{g}_{\alpha}) \geq 2$ or that there exists some $\gamma \neq \alpha$ also in $S$.  

Suppose $\mbox{dim}(\mathfrak{g}_{\alpha}) = 1$.  The assumption that $\mathfrak{g}$ has no $\mathfrak{sl}_2(\R)$-factor implies that there exists some nonzero root $\delta \neq \alpha$ such that $\delta(J) \neq 0$.  We may assume $\delta(J) < 0$.  If $\delta(H) < 0$, then $\delta \in S$, as desired.  So suppose that $\delta(H) \geq 0$.  Now let 
$$\gamma = - \sigma_{\alpha}(\delta) = - \delta - \delta(J) \alpha$$
Then 
$$\gamma(J) = - \delta(J) - \delta(J) \alpha(J) = - \delta(J)(1 - 2) = \delta(J) < 0$$
and 
$$\gamma(H) = - \delta(H) - \delta(J) \alpha(H) < 0$$
so $\gamma \in S$.
\end{proof}

\section{Nonproper semisimple actions: Proof of Theorem \ref{main.result}}
 
 {\it Reduction of the proof.} 
 As discussed in the previous section, it is proved in \cite{ADZ} that $G$ has an orbit $O$ with a noncompact stabilizer. It is also
 proved (Theorem 1.5) that if $O$ is Lorentzian, then the situation is exactly as described in Theorem \ref{main.result}.  If an orbit is not Lorentzian, then it is degenerate, fixed, or Riemannian.  Riemannian isometries always have compact isotropy, so the existence of noncompact stabilizers implies that either $O$ is degenerate, or some simple local factor fixes it pointwise.
The proof would be finished using the following two propositions, which state that in either situation, there is a nonproper Lorentzian orbit. 
 
\begin{proposition} \label{fixed}  Let $O$ be a $G$-orbit on which some simple noncompact normal subgroup $G_1$ acts trivially.  Then $G$ has (near $O$) Lorentzian orbits with noncompact isotropy.
\end{proposition}

\begin{proposition}  \label{degenerate} Let $O$ be a degenerate $G$-orbit with noncompact isotropy.  Then, $G$ has (near $O$) Lorentzian orbits with noncompact isotropy.
\end{proposition}

\begin{proof} (of Proposition \ref{fixed})  Let $x$ be a point of $O$.  The isotropy representation $\Phi : G_1 \rightarrow O(T_xM)$ is faithful and proper.  No noncompact simple subgroup of $O(T_xM) \cong O(1,d-1)$, can preserve an isotropic line.  Therefore, Theorem 1.1 of \cite{BZ} applies, so $\Phi(G_1)$ preserves an $(n+1)$-dimensional Lorentz subspace $V_1$ of $T_xM$ for some $2 \leq n \leq d-1$; further, $\Phi(G_1)$ contains a subgroup conjugate to the standard copy of $O^0(1,n)$ in $O(1,d-1)$.  Because $\Phi(G_1)$ is simple and $G$ has no local $SL_2(\R)$-factors, $G_1$ is locally isomorphic to $O(1,n)$, and $n \geq 3$.

Any spacelike vector $v \in V_1$ has Lorentzian $\Phi(G_1)$-orbit isometric to $dS_n$.  The stabilizer of $v$ is locally isomorphic to $O(1,n-1)$; in particular, it is noncompact.  The $G_1$-orbit of $exp_x(v)$ is again of Lorentzian type, by the Gauss lemma (see \cite{ON} 5.1). The $G$-orbit of $exp_x(v)$ is Lorentzian because it contains a Lorentzian submanifold, and this orbit has noncompact isotropy, as desired.
\end{proof}

The remainder of the section is devoted to the proof of Proposition \ref{degenerate}; the orbit $O$ will be assumed degenerate with noncompact isotropy below.  
 
\subsection{The asymptotic geodesic hypersurface ${\mathcal F}_x$}

\begin{fact} For  $x \in O$, let ${\R} n_x$ be the null direction in $T_xO$. Then, the orthogonal $n_x^\perp$ is tangent to a lightlike geodesic hypersurface ${\mathcal F}_x$ (defined in a neighborhood of $x$).
\end{fact}

\begin{proof}
Let $X$ be a nilpotent element of $\mathfrak{g}(x)$ given by Proposition \ref{propJ}, and consider the isometry $f = e^{tX}$, for some $t \neq 0$.  The derivative $D_xf(n_x) = n_x$.  Recall that $d$ is the dimension of $M$; let $g$ the Lorentz metric.  The graph $Graph(f) \subset M \times M$ is an isotropic totally geodesic $d$-dimensional submanifold of $M \times M$, equipped with the metric $g \oplus (-g)$. The graphs $Graph(f^m)$ converge to $E$, a $d$-dimensional, isotropic, totally geodesic submanifold, which is no longer a graph, since $f^m$ is divergent (see \cite{zecc} or \cite{DAG} 7.4). The intersection $E \cap (\{x\} \times M)$ is nontrivial, but has dimension at most $1$, because it is isotropic, and $M$ is Lorentzian.  Therefore, the projection  ${\mathcal F}_x$ of $E$ is a totally geodesic hypersurface in $M \times \{x\}$.  Because the derivative $D_x(f^m)$ fixes $n_x$ for all $m$, the vector $(n_x,n_x) \in T_{(x,x)}E$.  Because $T_{(x,x)}E$ is isotropic, its intersection with ${\bf 0} \times T_xM$ is exactly $\R n_x$.  Then $(n_x, {\bf 0}) \in T_{(x,x)} E$, so the projection $T_x (\mathcal{F}_x) = n_x^{\perp}$, as desired. 
\end{proof}

\begin{fact} The hypersurface ${\mathcal F}_x$ carries a 1-dimensional 
foliation ${\mathcal C}$ such that: 

\begin{enumerate}
\item{Any isotropic curve in ${\mathcal F}_x$ is tangent to a leaf of ${\mathcal C}$.}
\item{Each leaf of $\mathcal{C}$ is an isotropic geodesic.}
\item{The (local) quotient space ${\mathcal F}_x / {\mathcal C}$ inherits a Riemannian 
metric, infinitesimally preserved by the  elements of $\mathfrak{g}$ preserving ${\mathcal F}_x$.}
\end{enumerate}
\end{fact}

\begin{proof}  At any point $y$ of a degenerate hypersurface $\mathcal{F}$, there exists a unique tangent isotropic direction $C_y$.  These lines determine a characteristic 1-dimensional foliation ${\mathcal C}$ of $\mathcal{F}$, proving (1).  Since $\mathcal{F}_x$ is totally geodesic, (2) follows.  For (3), it is known  (see, for instance, \cite{Ze2}) that  ${\mathcal C}$ is transversally Riemannian if and only if $\mathcal{F}_x$ is totally geodesic. Here, transversally Riemannian means that the flow along any parameterization of ${\mathcal C}$ preserves the induced degenerate metric, or equivalently, that the degenerate metric can be projected as a Riemannian metric on the (local) quotient space $\mathcal{F}_x/{\mathcal C}$.  Finally, if an isometric flow locally preserves $\mathcal{F}_x$, then it induces a local diffeomorphism of $\mathcal{F}_x / \mathcal{C}$ that is obviously an isometry by construction of the Riemannian metric.
\end{proof}

\begin{fact}
The subalgebra $\mathfrak{s}(x)$ preserves the isotropic geodesic ${\mathcal C}_x$, so it preserves $\mathcal{F}_x$.
\end{fact}

\begin{proof} Indeed, any $Y \in \mathfrak{s}(x)$ has $\overline{Y}(x)$ isotropic, and hence the whole $Y$-orbit of $x$
is isotropic. But, as stated above, isotropic curves of ${\mathcal F}_x$ are contained in leaves of ${\mathcal C}$---that is, all $Y$-orbits through $x$ are contained in ${\mathcal C}_x$.  The image of $\mathcal{C}_x$ by any element of the one-parameter group generated by $Y$ is an isotropic geodesic tangent to $\mathcal{C}_x$ at some point, thus equals $\mathcal{C}_x$. 
\end{proof}

\begin{fact}  The action of $\Sigma_{\alpha \in S} \mathfrak{g}_\alpha$ on the tangent space of the (local) quotient space ${\mathcal F}_x / {\mathcal C}$ at the point corresponding to ${\mathcal C}_x$ is trivial.
\end{fact}

\begin{proof}
The tangent space to $\mathcal{F}_x / \mathcal{C}$ at $\mathcal{C}_x$ is identified with $n_x^\perp / {\R}n_x$.   
Note that the subspace $\R J + \Sigma_{\alpha \in S} \mathfrak{g}_{\alpha}$ as in Proposition \ref{propJ} is in fact a subalgebra of $\mathfrak{s}(x)$.
We have a representation $\rho$ of $\R J + \Sigma_{\alpha \in S} \mathfrak{g}_\alpha$ into the orthogonal algebra of $n_x^\perp / \R n_x$, which is endowed with a positive definite inner product. But
in such an orthogonal algebra, an equality $[\rho(J), \rho(Y)] = \lambda \rho(Y)$, becomes trivial---that is $\rho(Y) = 0$ (since $\lambda \neq 0$);
\end{proof}

\begin{corollary} \label{trivial.quotient}  $\Sigma_{\alpha \in S} \mathfrak{g}_{\alpha} $ acts trivially on the (local) 
quotient space ${\mathcal F}_x /{\mathcal C}$. That is,  $\Sigma_{\alpha \in S} \mathfrak {g}_{\alpha} $
preserves individually each leaf of ${\mathcal C}$.
\end{corollary}

\begin{proof} The action of $\Sigma_{\alpha \in S} \mathfrak{g}_{\alpha}$ on ${\mathcal F}_x/{\mathcal C}$
is trivial, since it is a Riemannian action with a  fixed a point and a trivial derivative at it.
\end{proof}

\begin{corollary} 
\label {any.point} 
Any point of ${\mathcal F}_x$ has a noncompact isotropy algebra. 
\end{corollary}

\begin{proof}  Indeed, $\Sigma_{\alpha \in S} \mathfrak{g}_{\alpha} $ has 
dimension $\geq 2  $ and has orbits of dimension 1. Therefore, stabilizers are nontrivial.  They are not compact since all elements of $\Sigma_{\alpha \in S} \mathfrak{g}_{\alpha}$ are nilpotent.
\end{proof}

\begin{fact} \label{fixed.points}  Let $\Gamma $ be the set of fixed points of  $\Sigma_{\alpha \in S} \mathfrak{g}_{\alpha}$ in ${\mathcal F}_x$. Then, $\Gamma$ has an empty interior (in ${\mathcal  F}_x$). In particular, the orbit of any point of ${\mathcal F}_x  - \Gamma$ under $\Sigma_{\alpha \in S} \mathfrak{g}_{\alpha}$ locally coincides with its ${\mathcal C}$-leaf. 
 \end{fact}

\begin{proof} No element of $\Sigma_{\alpha \in S} \mathfrak{g}_{\alpha}$ 
can fix points of an open subset of ${\mathcal F}_x$.  Indeed, in general, a Lorentz transformation fixing one point 
and acting trivially on a tangent lightlike hyperplane at that point has a trivial derivative, and is therefore trivial. 
\end{proof}

\begin{corollary} \label{non.spacelike}  No point of ${\mathcal F}_x$ has a spacelike $G$-orbit.
\end{corollary}

\begin{proof}  If a point $y \in {\mathcal F}_x$ has a spacelike orbit, then all orbits of points in a neighborhood of $y$ are spacelike, as well.  However, any neighborhood of $y$ meets ${\mathcal F}_x - \Gamma$; orbits of points in here cannot be spacelike, because they contain at least one isotropic geodesic. 
\end{proof}

\subsection{End of the proof of Proposition \ref{degenerate}}

\begin{fact} \label{deg.not.contained} The degenerate orbit $O$ cannot be locally contained in ${\mathcal F}_x$---that is, $\mathcal{F}_x \cap O$ does not contain an open subset of $O$.
\end{fact}

\begin{proof} Suppose $O$ is locally contained in $\mathcal{F}_x$.  Then the group $G$ locally preserves ${\mathcal F}_x$.  From Corollary  \ref{trivial.quotient}, the infinitesimal action of $\mathfrak{g}$ on the quotient space $Q = {\mathcal F}_x /{\mathcal C}$ is not faithful.  More precisely, any factor $\mathfrak{b}$ of $\mathfrak{g}$ which contains an element like $J$ (in Proposition \ref{propJ}) must act trivially on $Q$. However, orbits of $\mathfrak{b}$ cannot be 1-dimensional, since $\mathfrak{g}$ has no $\mathfrak{sl}_2(\R)$-factor. Therefore, $\mathfrak{b}$
acts trivially on ${\mathcal F}_x$.  As in the proof of Fact \ref{fixed.points}, this implies $\mathfrak{b}$ acts trivially on $M$.
\end{proof}

Now, from Corollary  \ref{any.point}, the proof of Proposition \ref{degenerate} would be finished once one proves that there is a point of ${\mathcal F}_x$ with a Lorentz orbit. It suffices
to show existence of nondegenerate orbits, since from Corollary \ref{non.spacelike}, points of ${\mathcal F}_x$
cannot have spacelike orbits. Assume, for a contradiction, that all $G$-orbits of points of ${\mathcal F}_x$ are degenerate. For any $y \in \mathcal{F}_x - \Gamma$, the orbit $Gy$ locally contains the isotropic geodesic $\mathcal{C}_y$ by Fact \ref{fixed.points}; therefore, $T_y(Gy) \subset n_y^{\perp} = T_y({\mathcal F}_x)$.  In other words, for any Killing field $X \in \mathfrak{g}$ and any $y \in \mathcal{F}_x - \Gamma$, the evaluation $X(y)$ is tangent to $\mathcal{F}_x$, so $X$ defines a vector field in this open subset of $\mathcal{F}_x$.  The flow of any $y \in \mathcal{F}_x - \Gamma$ along $X$ for sufficiently short time is again contained in $\mathcal{F}_x$.  In particular, since $x \in \mathcal{F}_x - \Gamma$, the orbit $O$ is locally contained in $\mathcal{F}_x$, contradicting the previous fact.

 %%  Let $X$ be a nilpotent element of $\mathfrak{g}(y)$.  As in the proof of Proposition \ref{propJ}, there is an $\R$-split element $A$ such that $A \in \mathfrak{s}(y)$ and $[J,X] = -2X$.   Then Corollary \ref{trivial.quotient} applies to give that $X$ is trivial on $\mathcal{F}_y / \mathcal{C}$.  As in Lemma \ref{3props}, the isotropy representation of $\mathfrak{g}(y)$ is trivial on the line $T_y \mathcal{C}_y$.  We conclude, as in Fact \ref{deg.not.contained}, that $X$ is trivial, a contradiction.  This completes the proof of Proposition \ref{degenerate}.

\begin{corollary} (from proof) There is a simple local factor $G_1$ of $G$ for which the $G_1$-orbit of $x$ is a point or degenerate with 
noncompact stabilizer. In other words, if $G$ has a nonproper orbit that is either a point or degenerate, then 
a subgroup $G_1$ locally isomorphic to $O(1,n)$ or $O(2,n)$ has an orbit with the same properties. 
\end{corollary}

\begin{proof}  
First suppose that $G$ has a nonproper fixed or Riemannian orbit $O_1$.  Then any local factor $G_1$ in the kernel of the restriction to $O_1$ has the desired properties.

Now suppose that $G$ has a nonproper degenerate orbit.  We have seen that some nilpotent elements stabilizing a point in the degenerate orbit $O$ stabilize a point $y$ with Lorentzian orbit. But, from Theorem \ref{adz}, a Lorentz orbit can be nonproper only if there is a local factor $G_1$ acting nonproperly with $O_1= G_1y$ Lorentzian of constant curvature; moreover, there is a warped product $L \times_w O_1$ preserved by $G$, and $G$ splits up to finite cover as $G_2 \times G_1$.  If $X \in \mathfrak{g}(y)$, the projection of $X$ on $\mathfrak{g}_2$ generates a precompact 1-parameter group, so $X$ cannot be nilpotent unless it belongs to $\mathfrak{g}_1$.  We infer from this that $G_1$ acts nonproperly on the degenerate orbit $O = Gx$.  Because $G_1$ is simple and the stabilizer $G_1(x)$ is noncompact, the orbit $G_1x$ must be degenerate or one point.  

\end{proof}

\section{Degenerate stabilizers: Proof of Theorem \ref{degstabs} }

Now as in Theorem \ref{degstabs}, assume $\mathfrak{g} \cong \mathfrak{o}(1,n)$ or $\mathfrak{o}(2,n)$ for some $n \geq 3$ and that $\mathfrak{g}$ has a degenerate orbit on $M$ with noncompact stabilizer $\mathfrak{g}(x)$.  It is proved in the first section below that $\mathfrak{g} \cong \mathfrak{o}(1,n)$ and in the second that $\mathfrak{g}(x) \cong \mathfrak{o}(n-1) \ltimes \R^{n-1}$.  The final point of the theorem follows from these two.

\subsection{Excluding $\mathfrak{o}(2,n)$}
\

Let $J, \mathfrak{a},$ and $S \subset \Delta$ be as in Proposition \ref{propJ}.  Observe that there can be at most one negative root $\alpha$ with $\mathfrak{g}_{\alpha} \cap \mathfrak{g}(x) \neq {\bf 0}$.  For if $X \in \mathfrak{g}_{\alpha} \cap \mathfrak{g}(x)$, then let $Y \in \mathfrak{g}_{- \alpha}$ be as in \cite{sam} 2.4.B, so $[X,Y] = H_{\alpha} \in \mathfrak{a}$.  Since $\overline{\mbox{ad}}(X)$ is nilpotent but skew-symmetric on $\mathfrak{g} / \mathfrak{s}(x)$ (Lemma \ref{3props} (3)), the $\R$-split element $H_{\alpha}$ must belong to $\mathfrak{s}(x)$.  Any $X^\prime \in \mathfrak{g}_{\beta}$ would give rise to $H_{\beta} \in \mathfrak{a} \cap \mathfrak{s}(x)$.  If $\alpha \neq \beta$, then for some $c \in \R$, the difference $H_{\alpha} - c H_{\beta}$ would be a nonzero $\R$-split element of $\mathfrak{g}(x)$, contradicting Lemma \ref{3props} (2).

Suppose that $\mathfrak{g} \cong \mathfrak{o}(2,n)$.  Let $\beta$ and $\gamma$ be distinct negative roots, each with $(n-2)$-dimensional root spaces.  The other negative roots are $\beta - \gamma$ and $\beta + \gamma$, with one-dimensional root spaces.  

First suppose $X \in \mathfrak{g}_{\beta} \cap \mathfrak{g}(x)$.  Let $L$ be a generator of $\mathfrak{g}_{-\beta - \gamma}$.  The adjoint $(\mbox{ad} X)^2(L)$ is a nonzero element of $\mathfrak{g}_{\beta - \gamma}$.  On the other hand, Lemma 3.2 (3) implies that any nilpotent element of $\mathfrak{g}(x)$ has nilpotence order $2$ on $\mathfrak{g} / \mathfrak{g}(x)$.  Then we would have $\mathfrak{g}_{\beta - \gamma} \cap \mathfrak{g}(x) \neq {\bf 0}$, a contradiction.

%% \begin{eqnarray*}
%% \left[ X,L \right] & = & W  \qquad \mbox{where} \ {\bf 0} \neq W \in \mathfrak{g}_{-\gamma}  \\
%% \left[ X,W \right] & = & S   \qquad \mbox{where} \ {\bf 0} \neq S \in \mathfrak{g}_{\beta - \gamma}  \\
%% \left[ X,S \right] & = & {\bf 0}
%% \end{eqnarray*}

%% The nilpotent subalgebra generated by $X$ and $L$ cannot be contained in $\mathfrak{s}(x)$ because $[ \mathfrak{g}(x), \mathfrak{s}(x) ] \subseteq \mathfrak{g}(x)$, and any nilpotent subalgebra of $\mathfrak{g}(x)$ is abelian. Let $\overline{L}$ be the image of $L$ modulo $\mathfrak{s}(x)$.  Because the adjoint of $X$ on $\mathfrak{g} / \mathfrak{s}(x)$ is skew-symmetric, $W$ must be in $\mathfrak{s}(x)$.  Because $X$ and $W$ generate a non-abelian nilpotent algebra, $W \notin \mathfrak{g}(x)$.  Then $cW - J \in \mathfrak{g}(x)$ for some $c \in \R$.  But $\overline{L}$ would be an eigenvector with real non-zero eigenvalue for $\ad(cW - J)$, contradicting \ref{3props} (1).

Therefore, $X$ cannot be in $\mathfrak{g}_{\beta}$.  The same argument shows $X$ cannot be in $\mathfrak{g}_{\gamma}$.  Proposition \ref{propJ} says that $\mathfrak{g}(x)$ contains a sum of negative root spaces with total dimension at least $2$, but since $\mathfrak{g}_{\beta - \gamma}$ and $\mathfrak{g}_{\beta + \gamma}$ are $1$-dimensional, and because only one can be contained in $\mathfrak{g}(x)$, we have a contradiction.

%% In fact, $\mathfrak{g}(x) \cap \mathfrak{g}_{\omega}$ must be ${\bf 0}$ for $\omega =  \pm \beta, \pm \gamma$.  

%% Next suppose $\alpha = \beta - \gamma$ and $X \in \mathfrak{g}_{\beta - \gamma} \subset \mathfrak{g}(x)$. The bracket $[X,\mathfrak{g}_{\gamma}] = \mathfrak{g}_{\beta}$, so the skew-symmetry condition \ref{3props} (3) forces $\mathfrak{g}_{\beta} \subset \mathfrak{s}(x)$.  From $[X,\mathfrak{g}_{-\beta}] = \mathfrak{g}_{-\gamma}$, also $\mathfrak{g}_{-\gamma} \subset \mathfrak{s}(x)$.  There is some nonzero $B+C \in \mathfrak{g}(x)$, where $B \in \mathfrak{g}_{\beta}$ and $C \in \mathfrak{g}_{-\gamma}$.  Since neither $\mathfrak{g}_{\beta}$ nor $\mathfrak{g}_{-\gamma}$ can intersect $\mathfrak{g}(x)$, both $B$ and $C$ must be nonzero.  Then, for a generator $L \in \mathfrak{g}_{-\beta - \gamma}$, the element $\mbox{ad}(B+C)^2(L)$ would be a nonzero element of  .

%% Last, suppose $X \in \mathfrak{g}_{\beta + \gamma}$.  A similar contradiction results from $[X, \mathfrak{g}_{-\beta}] = \mathfrak{g}_{\gamma}$ and $[X, \mathfrak{g}_{-\gamma}] = \mathfrak{g}_{\beta}$.

\subsection{Full stabilizer}

Now $\mathfrak{g}$ must be $\mathfrak{o}(1,n)$.  Let $\alpha$ and $J$ be as above, so $\alpha(J) < 0$ and $\mathfrak{g}_{\alpha} \subseteq \mathfrak{g}(x)$.  Let $\mathfrak{m}$ be the maximal compact subalgebra of the centralizer of $J$ in $\mathfrak{g}$; it is isomorphic to $\mathfrak{o}(n-1)$.  

Suppose $Y \in \mathfrak{s}(x) \cap \mathfrak{g}_{- \alpha}$.  By Lemma \ref{3props} (3),
$$ J \in [Y,\mathfrak{g}_{\alpha}] \subset [\mathfrak{s}(x),\mathfrak{s}(x)] \subset \mathfrak{g}(x)$$

But this contradicts Lemma \ref{3props} (2).  Therefore, $\mathfrak{s}(x) \cap \mathfrak{g}_{- \alpha} = {\bf 0}$.

On the other hand, since $\ad (X)$ is nilpotent for all $X \in \mathfrak{g}_{\alpha}$, Lemma \ref{3props} (3) forces 
$$\mathfrak{m} \subset [\mathfrak{g}_{- \alpha}, \mathfrak{g}_{\alpha}] \subset \mbox{ad}(\mathfrak{g}_{\alpha})(\mathfrak{g}) \subset \mathfrak{s}(x)$$

Since $\mathfrak{g} = \mathfrak{g}_{- \alpha} + \mathfrak{m} + \R J + \mathfrak{g}_{\alpha}$, the algebra $\mathfrak{s}(x)$ is exactly $\mathfrak{m} + \R J + \mathfrak{g}_{\alpha}$.  Suppose there were $X = cJ + M \in (\R J + \mathfrak{m}) \cap \mathfrak{g}(x)$ for some nonzero $c \in \R$. The subspace $\mathfrak{g}_{-\alpha}$ is $\mbox{ad}(X)$-invariant and maps onto $\mathfrak{g}/\mathfrak{s}(x)$.  But $\mbox{ad}(X)$ is clearly not skew-symmetric here, contradicting Lemma \ref{3props} (3).  Therefore, $\mathfrak{g}(x)$ is exactly $\mathfrak{m} + \mathfrak{g}_{\alpha}$, which is isomorphic as a Lie algebra to $\mathfrak{o}(n-1) \ltimes \R^{n-1}$.

\section{Global AdS warped product: Proof of Theorem \ref{warped} }
\label{section.global.ads.warped}

Suppose $G$ is locally isomorphic to $O(2, n)$, $n \geq 3$, with finite center, and $G$ acts isometrically on $M$.  By the argument of Kowalsky and
the assumption of no local $SL_2(\R)$-factors, we know that there is a $G$-orbit with
noncompact stabilizer (see Section \ref{section.roots.in.isotropy}).  By Theorem \ref{degstabs}, any $G$-orbit with noncompact stabilizer is Lorentzian.  Then by Theorem \ref{adz}, a neighborhood of some $G$-orbit is a warped product of the form $L \times_w S$, where $S$ is isometric to $AdS_{n+1}$ up to finite cover.  The set of orbits having a neighborhood isometric to $L \times_w S$, for some Riemannian manifold $L$ and $w : L \rightarrow \R^+$, is open. Let us prove that this set is also closed, and thus equals the whole of $M$.  A limit $O$ of a sequence $O_k$ of such orbits is a non-Riemannian orbit $O$ of dimension $\leq n + 1$.  Suppose that such a limit $O$ has compact isotropy $G(x)$ for $x \in O$.  Then $G(x)$ is contained in a maximal compact subgroup $K$ of $G$.  The Lie algebra $\mathfrak{k} \cong \mathfrak{o}(2) \times \mathfrak{o}(n)$, which has codimension $2n$ in $\mathfrak{g}$.  Since $n \geq 3$, this is impossible.  Therefore, for any $x \in O$, the stabilizer $G(x)$ is noncompact.  From Theorem \ref{degstabs}, $O$ cannot be degenerate; hence, it is Lorentzian.  Then by the \cite{ADZ} result (Theorem \ref{adz} above), a neighborhood of this orbit is isometric to a warped product.  Any orbit of $M$ has a neighborhood isometric to $L \times_w S$, for some $L$ and $w$. 

From this, one sees in particular that the $G$-action determines a foliation ${\mathcal O}$.  In addition, ${\mathcal O}$ admits an orthogonal foliation ${\mathcal O}^\perp$. We will use the $G$-action to show that the pair of  foliations ${\mathcal O}$ and ${\mathcal O}^\perp$ arise from a global warped product of the form $L \times_w AdS_{n+1}$ on a finite cover of $M$. 

Choose a point $x_0 \in M$.  Let $O_0$ and $O_0^\perp$ be the leaves of $x_0$ in the foliations $\mathcal{O}$ and $\mathcal{O}^{\perp}$, respectively.  Let $H_0$ be the stabilizer of $x_0$.  Note that $O_0^{\perp}$ is a component of the fixed set $Fix(H_0)$.  The full $Fix(H_0)$ is $N(H_0) O_0^{\perp}$, where $N(H_0)$ is the normalizer of $H_0$ in $G$.  It is well known that $N(H_0)/H_0$ is finite; then since $O_0^{\perp}$ is $H_0$-invariant, $Fix(H_0)$ has finitely-many components, each isometric to $O_0^{\perp}$.

Let $i$ and $i^\perp$ denote the respective inclusions of $O_0$ and $Fix(H_0)$ in $M$.  Let $G$ act on $Fix(H_0) \times O_0$
by $g(x, y) = (x, gy)$.  Define a mapping $\phi: Fix(H_0) \times O_0 \to M$, by $\phi(x, gx_0) = g(i^\perp(x))$. One sees that $\phi$  is well-defined, and it is in fact the $G$-equivariant extension of the inclusions. 

Next, $\phi$ is a covering map.  Clearly $\phi$ is a local diffeomorphism.  It is also easy to see that $\phi$ is surjective: the orbit of any $y \in M$ is homothetic to $S$.  Let $H_y$ be its stabilizer.  There is some $g \in G$ conjugating $H_y$ to $H_0$.  Then $gy \in Fix(H_0)$, and $y = \phi(gy,g^{-1}x_0)$.  Finally, $\phi$ is everywhere $N$-to-$1$, where $N = |N(H_0)/H_0|$, because every $G$-orbit in $M$ is homothetic to $S$.  An $N$-to-$1$ surjective local diffeomorphism is a covering map.

The submanifold $Fix(H_0)$ is Riemannian.  Let $L = Fix(H_0)$ with the metric pulled back by $i^{\perp}$.  By $G$-equivariance of $\phi$, all leaves $\phi(L \times \{ y \})$ are isometric.  Also by $G$-equivariance, the metrics along $G$-orbits $\phi(\{ y \} \times S)$ are all homothetic, with homothety factor depending only on $y \in L$.  Since $S$ and $AdS_{n+1}$ have a common finite cover, $M$ is isometric, up to finite cover, to a warped product $L \times_w AdS_{n+1}$. $\Diamond$

%% For this, one proves that $\phi$ satisfies the path lifting 
%% property, that is, if $c: t \in [0, 1]  \to M$ is a smooth  path, and 
%% $\tilde{x_0}$ such that $\phi(\tilde{x_0}) = c(0)$,  then, there exists 
%% a smooth path $\tilde{c}$, such that $\phi \circ \tilde{c} = c$, and
%% $\tilde{c}(0) = \tilde{x_0}$. For this goal, one argues by contradiction: consider 
%% $t_0$ the first $t$ when the lift becomes impossible. One gets a contradiction
%% if $t_0 <1$, by seeing to the product topology near $c(t_0)$, and using
%% that the $G$-action sends diffeomorphically a leaf of 
%% ${\mathcal O}^\perp$ to a leaf of ${\mathcal O}^\perp$.

%% Therefore, $M$ is a quotient of $Fix(H_0) \times O_0$ by a discrete isometry group $\Gamma \cong N(H_0)/H_0$. But, since $G$
%% acts on $M$, the projection of $\Gamma$ on $Isom (O_0)$ commutes with 
%% the $G$ action. But, it is known that the $O(2, n)$-action on $AdS_n$ has a trivial 
%% centrallizer, and hence $\Gamma$ is contained in $Isom (O_0^\perp)$.  Therefore
%% $M$ is a global product $(O_0^\perp/\Gamma) \times O_0 $. Actually, 
%% this implies that the $O^\perp$-leaf of $x_0$ is in fact
%% $ O_0^\perp/ \Gamma$, that is $\Gamma$ is trivial. Finally, by construction, this
%% product  is orthogonal, and it is easily seen to be a  warped product. $\Diamond$

\section{Full description: Proof of Theorem \ref{full} }
 
Item (1) of Theorem \ref{full} follows from Theorems \ref{main.result} and \ref{warped}.  

In this section, we consider the case in which $G_1 \cong O(1,n)$.  The first point of item (2) follows from Theorem \ref{main.result}.  In the first subsection below, we will deduce the second point of (2) from Theorem \ref{degstabs}.  Next, we address the last point of item (2), to obtain a warped product near Riemannian $G_1$-orbits, as well as the global decomposition $M = U \sqcup \partial U \sqcup V$.  These two subsections will complete the proof of Theorem \ref{full}.

%% Last, we will show that all factors other than $G_1$ act properly on $M$.  

\subsection{Degenerate orbits}

 In this subsection, we assume $G$ is locally isomorphic to $O(1,n)$ for some $n \geq 3$, and that $G$ acts nonproperly on $M$.
Let $U$ be the set of points having $G$-orbit homothetic, up to finite cover, to $dS_n$.  
In order to deduce the second point of (2) from Theorem \ref{degstabs}, it suffices to prove the following lemma:

\begin{lemma}
\label{degenerate.limits}
Let $x \in \partial U$.  Then the $G$-orbit of $x$ is either a fixed point or degenerate with noncompact isotropy.
\end{lemma}
  
\begin{proof}
The $G$-orbit of $x$ cannot be Riemannian and has dimension at most $n$.  Suppose $G(x)$ is compact, so it is contained in a maximal compact subgroup $K$ of $G$.  The Lie algebra $\mathfrak{k} \cong \mathfrak{o}(n)$ and has codimension $n$, so $\mathfrak{g}(x) \cong \mathfrak{o}(n)$; further, $Gx$ is either Lorentzian or degenerate and $n$-dimensional.  In either case, the isotropy representation of $\mathfrak{g}(x)$ has a $1$-dimensional invariant subspace tangent to $Gx$ and an $(n-1)$-dimensional complementary representation in $T_x(Gx)$, which is necessarily trivial.  

On the other hand, $x$ is the limit of a sequence $x_i \in U$ for which the isotropy $\mathfrak{g}(x_i)$ is trivial on $T(Gx_i)^\perp$, because $\mathfrak{o}(1,n-1) \cong \mathfrak{g}(x_i)$ has no nontrivial representation in $\mathfrak{o}(d-n)$.  The limit $\lim \mathfrak{g}(x_i) \subseteq \mathfrak{g}(x)$, and, because these subalgebras have the same dimension, they must be equal.  Then continuity of the action implies that $\mathfrak{g}(x)$ is trivial on $T(Gx)^{\perp}$.  Now, whether $Gx$ is Lorentzian or degenerate, the isotropy $\mathfrak{g}(x)$ is trivial on all of $T_xM$.  But isometries fixing a point and having trivial derivative at that point are trivial, so we have a contradiction.

Now $G(x)$ is noncompact. If the orbit $Gx$ were Lorentzian, then the result of \cite{ADZ} (Theorem \ref{adz} above) would give that $x \in U$.  Therefore, the $G$-orbit of $x$ is fixed or degenerate.
\end{proof}

\subsection{Riemannian orbits}

 As above, we assume in this subsection that $G$ is locally isomorphic to $O(1,n)$ for some $n \geq 3$, and that $G$ acts nonproperly on $M$.  By Theorem \ref{main.result} (3), the set $U$ defined above is open and nonempty.  The following proposition implies the last point of item (2) in Theorem \ref{full}, as well as the claimed decomposition of $M$ into de Sitter orbits, fixed points and light cone orbits, and hyperbolic orbits.

See Chapter 6 of \cite{Kowthesis} for local versions of many of the results below.
 
  \begin{propn}  
\label{partition} 
Let $U$ and $V$ be the points of $M$ having orbits homothetic up to finite cover of the de Sitter space $dS_n$ and the hyperbolic space ${\bf H}^n$, respectively. Then:
 \begin{enumerate}
 \item{Each point of $V$ has a $G$-invariant neighborhood isometric to a warped product $L \times_w {\bf H}^n$, for some Lorentz manifold $L$. In particular, $V$ is open.}
 
 \item{$\partial U = \partial V$, and it consists of all fixed points and orbits locally homothetic to the Minkowski light cone.}
 
 \item{$M = U \sqcup \partial U \sqcup V$}
\end{enumerate}
 \end{propn}
 
\emph{Notation}:  There exists a neighborhood of the ${\bf 0}$-section of $TM$ on which the exponential map is defined and injective on each fiber.  We fix one such neighborhood and denote it by $\Omega$ below.

We first collect some lemmas for the proof.
  
 \subsubsection{Maximal subalgebras of $\mathfrak{o}(1, n)$}   

We begin with two facts about certain natural subalgebras of $\mathfrak{o}(1,n)$.  A subalgebra $\mathfrak{h} \subset \mathfrak{g}$ is maximal if it is not contained nontrivially in another subalgebra: if $\mathfrak{h} \subseteq \mathfrak{h}^\prime \subseteq \mathfrak{g}$, then $\mathfrak{h}^\prime = \mathfrak{h}$ or $\mathfrak{h}^\prime = \mathfrak{g}$.

Consider the infinitesimal action of $\mathfrak{g} = \mathfrak{o}(1,n)$ on the projectivization ${\bf P}(\R^{1,n})$ of the standard representation with the standard inner product of type $(1,n)$.  For ${\bf v} \in \R^{1,n}$, denote by $\overline{{\bf v}} $ the image in ${\bf P}(\R^{1,n})$.  

\begin{itemize}
\item{If ${\bf v}$ is spacelike, then $\mathfrak{g}({\bf v}) = \mathfrak{g}(\overline{{\bf v}})$ is conjugate to $\mathfrak{g}({\bf e}_n)$, which is isomorphic to $\mathfrak{o}(1,n-1)$, with an obvious embedding in $\mathfrak{o}(1,n)$.}

\item{ For ${\bf v}$ timelike, $\mathfrak{g}({\bf v}) = \mathfrak{g}(\overline{{\bf v}})$ is conjugate to $\mathfrak{g}({\bf e}_1)$, which is isomorphic to $\mathfrak{o}(n)$, with an obvious embedding in $\mathfrak{o}(1,n)$. }

\item{ Finally, for ${\bf v}$ nonzero and isotropic, $\mathfrak{g}(\overline{{\bf v}})$ is conjugate to $\mathfrak{g}( \overline{{\bf e}_1 + {\bf e}_2} )$, which is isomorphic to the parabolic subalgebra $\mathfrak{p}$ as in Section \ref{section.isotropy.props}.  This subalgebra is isomorphic to $\mathfrak{sim}(n-1)$, the algebra of infinitesimal affine similarities of $\R^{n-1}$.  The annihilator of a nonzero isotropic ${\bf v}$ is a codimension-$1$ ideal of $\mathfrak{g}( \overline{{\bf v}})$, isomorphic to $\mathfrak{euc}(n-1)$, the algebra of infinitesimal affine isometries of $\R^{n-1}$. }

\end{itemize}

\medskip
\begin{lemma} \  
\label{maximal.subgroups} 

\begin{enumerate} 

\item{$\mathfrak{sim}(n-1)$ is a maximal subalgebra of $\mathfrak{o}(1,n)$.}

\item{$\mathfrak{euc}(n-1)$ is contained in exactly one maximal subalgebra, $\mathfrak{sim}(n-1)$.}
 
\item{Let $\mathfrak{h}$ be a proper subalgebra of $\mathfrak{o}(1,n)$ containing two different conjugates $\mathfrak{k}$ and $\mathfrak{k}^\prime$, each isomorphic to $\mathfrak{o}(n-1)$.  If $\mathfrak{k}$ and $\mathfrak{k}^\prime$ have no common fixed point in the light cone of $\R^{1,n}$, then $\mathfrak{h}$ is conjugate to $\mathfrak{o}(1,n-1)$ or $\mathfrak{o}(n)$.}

\end{enumerate} 
 \end{lemma}

 \begin{proof} \ 

\begin{enumerate}

\item{The subalgebra $\mathfrak{sim}(n-1)$ acts infinitesimally conformally on the projectivization of the light cone, which is conformally equivalent to $S^{n-1}$.  It is the infinitesimal stabilizer of one point and acts transitively on the complement of this point.  Then for any $X \notin \mathfrak{sim}(n-1)$, the algebra generated by $X$ and $\mathfrak{sim}(n-1)$ is transitive on $S^{n-1}$ and contains the full stabilizer subalgebra of each point.  It follows that any subalgebra properly containing $\mathfrak{sim}(n-1)$ is $\mathfrak{o}(1,n) \cong \mathfrak{conf}(S^{n-1})$, so $\mathfrak{sim}(n-1)$ is a maximal subalgebra.}

\item{Let $\mathfrak{h}$ be a maximal subalgebra of $\mathfrak{o}(1,n)$ containing $\mathfrak{euc}(n-1)$, and suppose that $\mathfrak{h}$ does not preserve any isotropic line in $\R^{1,n}$.  Because the corresponding connected group $H \subset O(1,n)$ is not compact, Theorem 1.1 of \cite{BZ} implies that $\mathfrak{h}$ preserves some Lorentz subspace of $\R^{1,n}$ and contains the full infinitesimal linear isometry algebra of this subspace. Because $\mathfrak{euc}(n-1)$ preserves no proper Lorentz subspace, $\mathfrak{h}$ must equal $\mathfrak{o}(1,n)$. }

\item{Let $H$ be the connected subgroup of $O(1,n)$ with Lie algebra $\mathfrak{h}$.  First suppose $H$ is compact, so it is conjugate to a subgroup of $O(n)$ containing two copies of $SO(n-1)$.  By an argument similar to that for $\mathfrak{sim}(n-1)$, the subalgebra $\mathfrak{o}(n-1)$ is maximal in $\mathfrak{o}(n)$, so $\mathfrak{h} = \mathfrak{o}(n)$. 
 
 Now suppose $H$ is not compact and does not preserve any isotropic line in $\R^{1,n}$.  Then by \cite{BZ} 1.1, $\mathfrak{h} \cong \mathfrak{o}(1,k) \times \mathfrak{l}$, for some $1 \leq k \leq n$, and some compact $L$.  Since $\mathfrak{o}(n-1) \subset \mathfrak{h}$, we must have $k \geq n-1$, and $\mathfrak{l}$ is trivial.  Because $\mathfrak{h}$ is a proper subalgebra, $\mathfrak{h} \cong \mathfrak{o}(1,n-1)$.

Finally, suppose that $\mathfrak{h}$ preserves an isotropic line, so it is conjugate to a subalgebra of $\mathfrak{sim}(n-1)$.  The subalgebras $\mathfrak{k}$ and $\mathfrak{k}'$ belong to the maximal traceless subalgebra $\mathfrak{euc}(n-1)$, in which case they fix a common null vector.}
\end{enumerate}
\end{proof}
 
\subsubsection{Fixed point sets} 

Recall that, for $H$ a subgroup of $G$, the fixed set is denoted $Fix(H)$.  Recall also that the exponential map is defined on $\Omega$ and injective on each fiber in it.  Each component of $Fix(H)$ is a totally geodesic submanifold of $M$.  Let $D$ be the points of $M$ having degenerate $G$-orbit with noncompact isotropy.  By Theorem \ref{degstabs}, any $x \in D$ has $\mathfrak{g}(x) \cong \mathfrak{o}(n-1) \ltimes \R^{n-1} \cong \mathfrak{euc}(n-1)$.
 
 \begin{lemma} \label{fixed.isotropy} Let $x \in D$, and denote by $O$ its $G$-orbit.  Near $x$, the fixed set $Fix(G(x))$ coincides with $\exp_x(T_xO^\perp \cap \Omega)$. It has dimension $d-n$, and intersects $O$ along the isotropic geodesic in $O$ through $x$.
 \end{lemma}  
 
\begin{proof}
 By Proposition \ref{degenerate}, we can approximate the degenerate orbit $O$ by Lorentzian orbits $O_i$ with noncompact isotropy.  By Theorem \ref{adz}, each $O_i$ is homothetic, up to finite cover, to $dS_n$.  Take $x_i \in O_i$ with $x_i \to x \in O$.  The limit $\lim \mathfrak{g}(x_i)$ is contained in $\mathfrak{g}(x) \cong \mathfrak{euc}(n-1)$; since $dim(\mathfrak{g}(x_i)) = dim(\mathfrak{g}(x))$ for all $i$, the limit equals $\mathfrak{g}(x)$.  
As in the proof of Lemma \ref{degenerate.limits}, triviality of $\mathfrak{g}(x_i)$ on $(T_{x_i}O_i)^\perp$ for all $i$ implies that, in the limit, the isotropy $\mathfrak{g}(x)$ is trivial on $(T_x O)^\perp$.  It is easy to see from the form of the isotropy representation of $\mathfrak{g}(x) \cong \mathfrak{euc}(n-1)$ that the maximal trivial subrepresentation in $T_xM$ is exactly $(T_xO)^\perp$.  Then the lemma follows with the exponential map and dimension counting. 
\end{proof}
 
 \begin{lemma}  \label{fixed.compact} Let $K$ be a maximal compact connected subgroup of $G(x)$. Then $Fix(K)$ has dimension $d- n +1$ and is of Lorentzian type. It meets $O$ along the isotropic geodesic in $O$ through $x$.
 \end{lemma}
 
 \begin{proof} 
The subgroup $K$ is conjugate to $SO(n-1) \subset Euc(n-1)$.  The isotropy representation of $K$ at $x$ fixes any nonzero normal vector $n_x \in (T_xO)^\perp \cap T_xO$ and acts irreducibly on an $(n-1)$-dimensional complementary spacelike subspace $L_x \subset T_xO$.  Then $Fix(K) \cap O = N_x$, where $N_x$ is the isotropic geodesic in $O$ through $x$. Now, $L_x \oplus (T_xO)^\perp$ has codimension $1$ in $T_xM$.  There is a $K$-invariant $1$-dimensional complementary representation, which is necessarily trivial.  We have therefore a $(d-n+1)$-dimensional Lorentzian subspace $L_x^\prime$ containing $(T_xO)^\perp$, complementary to $L_x$, on which $K$ acts trivially.  This subspace is exactly the maximal trivial subrepresentation of $K$ on $T_xM$.
\end{proof}

\begin{lemma}  \label{exp} Let
$$X_x = \cup\{ T_x(Fix (K))\setminus \{ {\bf 0} \} \ : \  K \subset G(x) \ \mbox{maximal compact connected} \}$$
 
Then $X_x$ is open in $T_xM$. More precisely, 
$$X_x = ((T_xO)^\perp \times P^c) \cup ((T_xO)^\perp \setminus \{ 0 \})$$

where $P$ is a hyperplane in a subspace $F$ complementary to $(T_xO)^\perp$.

%% More exactly, if  $E= (T_xO)^\perp$, then, $X_x$
%%  contains a (product) set of the form $E \times V$, where $V$ is a cone in a subspace
%%  $F$ supplementary to $E$, that is,  $V$ is invariant under (positive and negative)  
%% scalar multiplication. Furthermore, $V$ has non-empty interior in $F$, but does not contain 0.
 \end{lemma}
 
 \begin{proof} 
Fix one maximal connected compact $K_0 \subset G(x)$, locally isomorphic to $SO(n-1)$.  Let $F_0$ be a $K_0$-invariant degenerate complement to $(T_xO)^\perp$, and let ${\bf v}_0$ be a nonzero isotropic vector in $F_0$.  Then $(T_xO)^\perp + \R {\bf v}_0 = T_x (Fix(K_0))$.  Choose a generator ${\bf n}_x$ of $(T_xO)^\perp \cap T_xO$; note that ${\bf n}_x$ is fixed by $G^0(x)$.  Now any other maximal compact, connected subgroup of $G(x)$ equals $g K_0 g^{-1}$ for some $g \in G^0(x)$.  The fixed subspace $T_x (Fix (g K_0 g^{-1})) = (T_xO)^\perp + \R g {\bf v}_0$.  As $g$ ranges over $G^0(x)$, the projection of $g {\bf v}_0$ to $F_0$ ranges over all vectors ${\bf v}$ with $\langle {\bf v}, {\bf n}_x \rangle =$ $\langle {\bf v}_0, {\bf n}_x\rangle$.  Then the projection of $\R G^0(x) {\bf v}_0$ to $F_0$ ranges over all ${\bf v}$ with $\langle{\bf v}, {\bf n}_x\rangle \neq 0$.  Now the lemma is proved, with $F = F_0$ and $P = F_0 \cap {\bf n}_x^\perp$.
\end{proof}

\subsubsection{Hypersurface of degenerate orbits}  Let $O$ be as above, a degenerate orbit with isotropy $\mathfrak{euc}(n-1)$.  We will next show that any neighborhood of a point of $O$ meets $U$ or $V$. Then we will show that $O$ lies in a hypersurface of degenerate orbits with the same isotropy, and, finally, that this hypersurface locally separates $U$ from $V$.

Let $X_x$ be as in Lemma \ref{exp} above.  Consider $X(O) = \cup \{ X_x \cap \Omega \ : \  x \in O\} $. This is an  open $G$-invariant set in the restriction $TM|_O$.  

\begin{lemma} \label{image} The image $exp(X(O)) \subset U \cup V$---that is, for any $y \in  exp(X(O))$, the stabilizer $\mathfrak{g}(y) \cong \mathfrak{o}(n)$ or $\mathfrak{o}(1,n-1)$.
\end{lemma}

\begin{proof}

Denote by $\Phi$ the restriction of the exponential map $\Phi = \exp: X(O) \to M$.  Each $X_x$ is open in $T_xM$, so $\Phi$ has  maximal rank, and is in particular an open map. Its image is an open set containing $O$ in its closure. Any point $y$ in this image is a regular value; the inverse image $S = \Phi^{-1}(y) \subset X(O)$ has dimension $dim (X(O)) - dim(M) = d $.  Because $exp$ is injective on each fiber of $\Omega$, and hence of $X(O)$, the inverse image $S$ projects injectively under $\pi : TM \rightarrow M$ onto some open set $W$ in $O$.  For any $z \in W$, there exists a maximal compact $K_z \subset G(z)$ such that $y \in Fix(K_z)$.  Because $n \geq 3$, there exist $z, z^\prime \in W$ which do not lie on a common isotropic geodesic of $O$, so that $K_z$ and $K_{z^\prime}$ are contained in $G(y)$ but have no common fixed point in $O$.  Now apply Lemma \ref{maximal.subgroups} (3).  
 \end{proof}

 Let $x \in O$, and let $R$ be a neighborhood of $x$ in $Fix(G(x))$.  Let $\widehat{O}$ be the union of all the $G$-orbits of all points of $R$.

\begin{lemma} \label{neighbourhood} For $R$ sufficiently small, $\widehat{O}$ is a degenerate hypersurface in which every orbit is degenerate with isotropy $\mathfrak{euc}(n-1)$; further, $\widehat{O} \cup U \cup V$ is a neighborhood of $\widehat{O}$.
\end{lemma}

\begin{proof}
From Lemma \ref{fixed.isotropy}, $R$ equals $\exp ((T_xO)^\perp \cap \Omega)$ near $x$ and has dimension $d-n$.  Since $\mathfrak{g}(x) \cong \mathfrak{euc}(n-1)$ lies in a unique maximal subalgebra of $\mathfrak{o}(1,n)$, the orbit of any $z \in R$ is either a fixed point or degenerate.  Shrink $R$ so that the $G$-orbit of any $z \in R$ is degenerate; in this case, $\mathfrak{g}(z) \cong \mathfrak{euc}(n-1)$.  The tangent space $T_z\widehat{O} = T_z(Gz) + T_z(Fix(G(x)) = T_z(Gz) + (T_z(Gz))^\perp$, because $z$ is also as in Lemma \ref{fixed.isotropy}; hence, $\widehat{O}$ is a degenerate hypersurface.
 
Now consider $X(\widehat{O}) = \cup\{ X(Gz) \ :\  z \in R\} \subset T\widehat{O}$, and let $\Phi$ be the restriction of $exp$ to $X(\widehat{O})$.  From Lemma \ref{image}, the image of the restriction $\Phi_z$ to $X(Gz)$ is open and contained in $U \cup V$.  The union $\Phi_z(X(Gz)) \cup Gz$ is a neighborhood of $Gz$.  Now taking 
 $$ \bigcup_{z \in R} \Phi_z(X(Gz)) \cup \widehat{O}$$

we obtain a neighborhood of $\widehat{O}$ contained in $\widehat{O} \cup U \cup V$.
\end{proof}

\begin{lemma}
\label{different.sides}
Let $\Phi$ be the restriction of the exponential map to $X(\widehat{O})$, in a sufficiently small neighborhood of the ${\bf 0}$-section.  One component of $\Phi(X(\widehat{O})) \backslash \widehat{O}$ lies in $U$, and the other is in $V$.
\end{lemma} 

\begin{proof}
From the last lemma, any point $x \in \widehat{O}$ has a neighborhood in which $dim(\mathfrak{g}(x)) = dim(\mathfrak{g}(y))$ for all $y$ in the neighborhood; let $N$ be this common dimension.  In this neighborhood, the map $y \mapsto \mathfrak{g}(y) \in Gr^N(\mathfrak{g})$ is continuous.  There are smaller neighborhoods $Z \subset W$ of $x$ such that $W$ is a normal neighborhood of each of its points, and, whenever $z, z' \in Z$ are such that $\mathfrak{g}(z) = \mbox{Ad}(g)(\mathfrak{g}(z'))$, then $gz' \in W$.  (\emph{Note}: We take $g$ to be the element of $\mbox{Ad}^{-1}(\mbox{Ad}(g))$ closest to the identity.)

Suppose for a contradiction that points on opposite sides of $\widehat{O}$ have the same orbit type.  Let $Z_1, Z_2$ be the components of $Z \backslash \widehat{O}$, and $W_1, W_2$ the components of $W \backslash \widehat{O}$.  Take $z \in Z_1$ and $z' \in Z_2$.  By assumption, $\mathfrak{g}(z)$ and $\mathfrak{g}(z')$ are conjugate in $G$.  Let $g$ be a conjugating element with minimal distance to the identity.  Then $z, gz' \in W$ are connected by a unique geodesic, which necessarily passes through $\widehat{O}$.  The common stabilizer $\mathfrak{g}(z) = \mathfrak{g}(gz')$ fixes the geodesic pointwise, so it fixes a point of $\widehat{O}$.  Whether $\mathfrak{g}(z)$ is isomorphic to $\mathfrak{o}(n)$ or to $\mathfrak{o}(1,n-1)$, neither embeds in the stabilizer of a point of $\widehat{O}$.
\end{proof}

\subsubsection{Proof of Proposition \ref{partition}} 

Suppose that $U \neq M$, so $\partial U \neq \emptyset$.  By Lemma \ref{different.sides}, not only is $V$ nonempty, but the interior $int(V) \neq \emptyset$.  The following fact will be key for the rest of the proof.

\begin{fact}
For $x \in int(V)$ with orbit $O$, the isotropy representation of $\mathfrak{g}(x)$ is trivial on $(T_xO)^\perp$.
\end{fact}

\begin{proof}
 The stabilizer $\mathfrak{g}(x) \cong \mathfrak{o}(n)$ acts on the orthogonal space $(T_x O)^\perp$ via a homomorphism $\mathfrak{o}(n) \to \mathfrak{o}(1,d-n-1)$.  There is a neighborhood $W$ of $(x,{\bf 0})$ in $(TO)^\perp$ on which the exponential map is a diffeomorphism onto its image.  For any $y \in exp_x (W \cap T_xM)$, the stabilizer $\mathfrak{g}(y) \subseteq \mathfrak{g}(x)$.  If all orbits near $O$ have the same dimension, then $\mathfrak{g}(x) = \mathfrak{g}(y)$ for all $y \in exp(W \cap T_xM)$.  Then the orthogonal representation of $\mathfrak{g}(x)$ is trivial.
\end{proof}

Now we can show that $\partial(int(V)) = \partial U$.  Suppose $x_n \in int(V)$ and $x_n \rightarrow x$.  Denote by $O$ the orbit of $x$. As usual, we have $\lim \mathfrak{g}(x_n) \subseteq \mathfrak{g}(x)$, and $dim(O) \leq n$.  The orbit of $x$ is either Riemannian, degenerate, or fixed.  

First suppose $O$ is Riemannian.  Then $\mathfrak{g}(x)$ is compact, so it is conjugate to a subalgebra of $\mathfrak{o}(n)$.  The dimension restriction forces $\mathfrak{g}(x)$ to be conjugate to $\mathfrak{o}(n)$, so $x \in V$.  As in the proof of \ref{degenerate.limits}, the isotropy $\mathfrak{g}(x)$ is trivial on $(T_xO)^\perp$, so orbits near $x$ have isotropy containing $\mathfrak{o}(n)$.  The orbit $O$ has a neighborhood consisting of Riemannian $G$-orbits.  Because $\mathfrak{o}(n)$ is a maximal subalgebra, orbits near $x$ all have isotropy $\mathfrak{g}(x)$.  Then $x \in int(V)$, a contradiction.

If $O$ is degenerate, the same argument as in Lemma \ref{degenerate.limits} implies that $\mathfrak{g}(x)$ is noncompact.  Then Proposition \ref{degenerate} implies that $x \in \partial U$.  If $x$ is a fixed point, then Proposition \ref{fixed} implies that $x \in \partial U$.  Now Lemma \ref{degenerate.limits} gives item (2).

Now $U \cup \partial U \cup int(V)$ is closed.  It is open by Lemma \ref{neighbourhood}, so it equals $M$.  Therefore, $int(V) = V$, and $M = U \sqcup \partial U \sqcup V$, proving item (3).

It remains to prove item (1), the warped product structure in $V$.  As in the proof in \cite{ADZ} (\S 2.3) of warped products around Lorentz orbits and as in Section \ref{section.global.ads.warped}, the key properties are irreducibility of the isotropy along orbits and triviality of the isotropy orthogonal to orbits.  
 
Fix $x \in V = int(V)$ with $G$-orbit $O$, and let ${\mathcal L}_x = \exp_x ((T_xO)^\perp \cap \Omega)$. Since $G(x)$ is irreducible on $T_xO$ and trivial on $(T_xO)^\perp$, the leaf ${\mathcal L}_{x}$ is exactly the set of common fixed points of $G(x)$ near $x$. Let $L = Fix(G(x)) \cap V$.  Define $\phi : L \times O \rightarrow V$ by $\phi(y,gx) = gy$.  This map is a local diffeomorphism because, for any $y \in Fix(G(x)) \cap V$, the fixed set $L$ coincides near $y$ with $exp((T_y(Gy))^\perp \cap \Omega)$.  If $\phi(y,gx) = \phi(y',g'x)$, then $g^{-1}g' \in N_G(G(x))$, and $g^{-1}g'y' = y$.  The stabilizer $G(x)$ is a maximal compact subgroup of the semisimple connected group $G$, so $N_G(G(x)) = G(x)$ (see \cite{Eber} 1.13.14 (4)). It follows that $\phi$ is injective.  Any $y \in V$ has $G(y)$ compact and $dim G(y) = dim G(x)$; because $G(x)$ is connected, there exists $g \in G$ such that $g G(y) g^{-1} = G(x)$.  Then $\phi(gy,g^{-1}x) = y$, so $\phi$ is surjective.  In conclusion, $\phi$ is a $G$-equivariant diffeomorphism $L \times O \rightarrow V$.

%% It is therefore a geodesic submanifold of $M$. One can do the same for all points of $O$ to obtain a family of geodesic submanifolds ${\mathcal L}$.  In fact, because $\mathfrak{g}(y) = \mathfrak{g}(x)$, for $y \in {\mathcal L}_{x}$, one obtains the same family $\mathcal{L}$ along the orbit of $y$.  Then ${\mathcal L}$ is a $G$-invariant foliation orthogonal to the $G$-orbits near $x$ and having geodesic leaves. 
%%  Take $L = {\mathcal L}_{x}$ and $N = O$.  Using the $G$-action, we can extend the canonical injections of $L$ and $N$ into $M$ to a natural mapping $p: L \times N \to  M$.  

The orbit $O$ is assumed to be homothetic, up to finite cover, to ${\bf H}^n$.  Because ${\bf H}^n$ is simply connected, it has no finite covers.  Because any finite subgroup of $Isom({\bf H}^n)$ has a common fixed point, there are no smooth finite quotients of ${\bf H}^n$.  Therefore, the orbit $O$, and every other orbit of $V$, is globally homothetic to ${\bf H}^n$.  The metric on $V$ pulls back by $\phi$ to $L \times O$.  The transverse fibers of the product are orthogonal by construction.  All $L$-fibers are isometric because they are related by the action of $G$.  Also by equivariance, $O$-fibers in $L \times O$ correspond to $G$-orbits in $V$, so each is homothetic to ${\bf H}^n$.  In conclusion, there exists $w : L \rightarrow \R^+$ such that $V$ is isometric to $L \times_w {\bf H}^n$. $\Box$

%% To prove that the metric on $V$ splits as a warped product, it suffices to show that the pullback metric on $L \times O$ has the form $h \bigoplus w m$ for some positive function $w$ on $L$.  This can be done following the same lines as in \cite{ADZ}, \S 2.3.  It is even simpler in our present situation.  First, by construction, the metric on $L \times O$ splits orthogonally.  The fact that the metric $h_{(l,n)}$ along a fiber $L \times \{ n \}$ depends only on $l$ follows from $G$-equivariance of $\phi$.  On the other hand, the transfer from a level $\{l\} \times N$ to $\{l^\prime\} \times N$ is equivalent to going from an orbit $O_x$ to another one $O_y$.  This motion is homothetic, since the $G$-action on ${\bf H}^n$ preserves exactly one metric up to a multiplicative constant. 

%% The map $L \times_w N \rightarrow M$ is a local isometry onto a neighborhood of $O$ in $M$.  This map can be made a global isometry by shrinking $L$ slightly.  Note that the stabilizer $G(x)$, locally isomorphic to $O(n)$, has finite index in its normalizer $N_G(G(x))$.  Let $\gamma_1, \ldots, \gamma_k$ be representatives for the nontrivial cosets in $N_G(G(x)) / G(x)$.  Now let 
%% $$ L = \mathcal{L}_x \backslash (\cup_{i=1}^k \gamma_i \overline{\mathcal{L}_x} \cap \mathcal{L}_x)$$

%% The new $L$ is an open subset of $\mathcal{L}_x$ containing $x$.  It is now easy to verify that $L \times N$ maps diffeomorphically and $G$-equivariantly to a neighborhood of $O$ in $M$.

The referee kindly communicated to us the following lemma.  It synthesizes arguments that we used in many proofs, especially in the proof of Proposition \ref{partition} above. Actually, this proof above includes the proof of the lemma. 
 
\begin{lemma}(the referee, anonymous)
Let $G$ act by isometries of a pseudo-Riemannian manifold $M$ with no isotropic orbits.
Suppose that, for all $x \in M$, the identity component $G^0(x)$ is conjugate to $S$, for $S$ a fixed closed, connected subgroup of $G$.
Assume that the adjoint representation of ${\mathfrak s}$
on ${\mathfrak g}/{\mathfrak s}$ is absolutely irreducible.
Let $L$ be the set of $S$-fixed points in $M$.
Then the $G$-action lifts to an $N$-to-1 cover $M'$ of $M$, where $N = [N_G(S):S]$.
Further, $M'$ is $G$-equivariantly isometric to a warped product $L \times_w G/S$, where $G/S$ carries a nondegenerate $G$-invariant pseudo-Riemannian metric, and $L$ carries the metric induced from $M$.
\end{lemma}

 The hypothesis that all stabilizers are conjugate in $G$ means that orbits are all isomorphic as homogeneous $G$-spaces.  Together with the irreducibility hypotheses, one obtains that the set of points having a given stabilizer $S^\prime$ is transverse to the orbit foliation, and in fact equals $Fix(S^\prime)$.  This situation is exactly that of the proof above.  Not only in the proof of (1) in Proposition \ref{partition}, but also in Section \ref{section.global.ads.warped}, as part of the proof of Theorem \ref{warped}, we prove a case of this lemma.  Note that the hypothesis that all stabilizers are conjugate is strong, and, in each of our two cases, fulfilling the hypotheses of the lemma is intertwined with the rest of the proof.

  {\bf Acknowledgements:}  We thank the referee for his many thoughtful comments and helpful suggestions.  The first two authors thank the \'Ecole Normale Sup\'erieure de Lyon for their hospitality.  During their stay, they received assistance from the ACI \emph{Structures g\'eom\'etriqes et trous noirs}.

{\small
\medskip
\noindent
Mohamed Deffaf, 
Facult\'e des math\'ematiques, USTHB \\
BP 32 El'Alia, Bab Ezzouar, Alger, Algeria.

\noindent
Karin Melnick,  
Department of Mathematics, Yale University \\
PO Box 208283, New Haven, CT, USA. 

\noindent
Abdelghani Zeghib, 
CNRS, UMPA, \'Ecole Normale Sup\'erieure de Lyon \\
46, all\'ee d'Italie,
69364 Lyon cedex 07,  France 
 
mhdeffaf@gmail.com \\
karin.melnick@yale.edu, www.math.yale.edu/\~{}km494/\\
Zeghib@umpa.ens-lyon.fr, www.umpa.ens-lyon.fr/\~{}zeghib/
}

\end{document}